\documentclass[12pt]{article}
\usepackage{amsmath,amssymb,enumerate,bbm}
\newtheorem{theorem}{Theorem}
\newtheorem{definition}{Definition}
\newtheorem{lemma}{Lemma}
\newenvironment{proof}{
  \noindent\textbf{Proof}\ }{\hspace*{\fill}
  \begin{math}\Box\end{math}\medskip}
\newcommand{\tmop}[1]{\ensuremath{\operatorname{#1}}}
\newenvironment{enumeratenumeric}{\begin{enumerate}[1.]}{\end{enumerate}}
\newtheorem{varremark}{Remark}
\newenvironment{remark}{\begin{varremark}\em}{\em\end{varremark}}

\begin{document}

\title{Divisor graphs have arbitrary order and size} \author{Le Anh Vinh\\School of Mathematics\\University of New South Wales\\Sydney 2052 Australia} \date{\empty} \maketitle

\begin{abstract}
  A divisor graph $G$ is an ordered pair $( V, E )$ where $V \subset
  \mathbbm{Z}$ and for all $u \neq v \in V$, $u v \in E$ if and only if $u \mid
  v$ or $v \mid u$. A graph which is isomorphic to a divisor graph is also called
  a divisor graph. In this note, we will prove that for any $n \geqslant 1$ and $0
  \leqslant m \leqslant \binom{n}{2}$ then there exists a divisor graph of order $n$
  and size $m$. We also present a simple proof of the
  characterization of divisor graphs which is due to Chartran, Muntean,
  Saenpholpant and Zhang. 
\end{abstract}

\section{Introduction}
The notion of divisor graph was first introduced by Singh and Santhosh \cite{hai}. A divisor graph $G(V)$ is an ordered pair $( V, E )$ where $V \subset
\mathbbm{Z}$ and for all $u \neq v \in V$, $u v \in E$ if and only if $u \mid v$
or $v \mid u$. A graph which is isomorphic to a divisor graph is also called a
divisor graph. The main result of this note is the following theorem.

\begin{theorem}\label{main}
For any $n \geqslant 1$ and $0
  \leqslant m \leqslant \binom{n}{2}$ then there exists a divisor graph of order $n$
  and size $m$.
\end{theorem}

To prove Theorem \ref{main}, we need the following characterization of divisor graphs is due to
Chartran, Muntean, Saenpholpant and Zhang \cite{mot}.

\begin{theorem}\label{theorem 1}
  A graph $G$ is divisor graph if and only if there is an orientation $D$ of
  $G$ such that if $( x, y ), ( y, z )$ are edges of $D$ then so is $( x, z
  )$.
\end{theorem}

The proof of Theorem \ref{theorem 1} in \cite{mot} is by induction on the order of the graph. For the completeness of this note, we will present a simple (and direct) proof of this theorem in Section 3. From Theorem \ref{theorem 1}, we introduce the definition of a divisor digraph which will be usefull in the proof of Theorem \ref{main}.

\begin{definition}
A digraph $G$ is a divisor digraph if and only if $(x,y), (y,z)$ are edges of $G$ then so is $(x,z)$.
\end{definition}

It is clear that if $G$ is a divisor digraph then the graph obtained by ignoring the direction of edges of $G$ is a divisor graph. 

\section{Proof of Theorem \ref{main}}

Suppose that $G = ( V, E )$ is a graph with vertex set $V = \{ v_1, \ldots,
v_n \}$ and size $m$. The degree of a vertex $v_i$ is the number of edges of
$G$ incident with $v_i$. Let $d_1 \geqslant d_2 \geqslant \ldots \geqslant
d_n$ be the vertex degrees in non-increasing order and let $f_G = | \{ i : d_i
\geqslant i \} |$. Let $e_i$ be the number of vertices with degree at least
$i$, that is $e_i = | \{ j : d_j \geqslant i \} |$. Then we have $e_1 \geqslant e_2 \geqslant \ldots \geqslant e_n$.
Moreover, we also have  \begin{equation} \label{s} d_i = | \{ j: e_j \geqslant i \} |. \end{equation}
We have the following lemmas.

\begin{lemma}\label{lemma 1}
  Let $n-1 = d_1 \geqslant \ldots \geqslant d_n \geqslant 1$ be a sequence of natural numbers
  and let $e_i = | \{ j : d_j \geqslant i \} |$. Suppose that
  \begin{equation} \label{1}\sum_{i = 1}^n d_i = 2 m, \end{equation}
  and
  \begin{equation}\label{2} \sum_{i = 1}^t d_i = \sum_{i = 1}^t ( e_i - 1 ) \end{equation}
  for all $1 \leqslant t \leqslant f = | \{ i : d_i \geqslant i \} |$. Then
  there exists a divisor graph of order $n$ and size $m$.
\end{lemma}

\begin{proof}
  We construct a digraph $G$ with vertex set $v_1, \ldots, v_n$ as follows.
  For $1 \leqslant i \leqslant f$, then $( v_i, v_j )$ is an edge of $G$ for
  $i + 1 \leqslant j \leqslant d_i + 1$. We first show that $G$ is a graph of
  size $m$. It suffices to show that $\deg ( v_i ) = d_i$ for all $i$ (where $deg(v_i)$ is the number of incident edges of vertex $v_i$, regardless of their directions). We have
  three cases.
  \begin{enumeratenumeric}
    \item Suppose that $1 \leqslant i \leqslant f$. Then it is clear from the
    construction that deg$( v_i ) = d_i .$
    
    \item Suppose that $i = f + 1$. Then $d_i = d_{f + 1} \leqslant f$ (by the
    definition of $f$). From (\ref{2}), we have $d_j = e_j - 1$ for $1 \leqslant j
    \leqslant f$. If $d_{f + 1} < f$ then $e_f \leqslant f$, or
    \[ f \leqslant d_f = e_f - 1 < f, \]
    which is a contradiction. Hence $d_{f + 1} = f$. For $1 \leqslant j
    \leqslant f$ then $d_j + 1 \geqslant d_f + 1 \geqslant f + 1 = i$. So $(
    v_j, v_i )$ is an edge of $G$ and $\deg ( v_i ) = f$. Thus, $\deg ( v_i )
    = d_i$.
    
    \item Suppose that $i > f + 1$. Then $( v_j, v_i )$ is an edge of $G$ if
    and only if $1 \leqslant j \leqslant f$ and $j + 1 \leqslant i \leqslant
    d_j + 1$. For $j > f$ we have $d_j \leqslant d_{f + 1} < f + 1 \leqslant
    i$. This implies that
    \[ \deg ( v_i ) = | \{ 1 \leqslant j \leqslant d : d_j \geqslant i - 1 \}
       | = | \{ j : d_j \geqslant i - 1 \}| = e_{i - 1} \]
    for all $i > f + 1$. From (\ref{s}), we have $d_i = e_i - 1$ for $1 \leqslant i \leqslant f$. For $i > f
+ 1$ then $d_j \geqslant i - 1$ or $e_j \geqslant i$ only if $j \leqslant f$
(since $d_{f + 1} = f < i - 1$). Hence
\begin{align*}
  e_{i - 1} &= | \{ j : d_j \geqslant i - 1 \} |\\
            &= | \{ 1 \leqslant j \leqslant f : d_j \geqslant i - 1 \} | \\
  				  &= | \{ 1 \leqslant j \leqslant f : e_j - 1 \geqslant i - 1 \}|\\
  					&= | \{ j : e_j \geqslant i \} | = d_i . 
\end{align*}
Thus, deg$( v_i ) = e_{i - 1} = d_i$ for all $i > f + 1$.
  \end{enumeratenumeric}
  
Therefore, we have deg$(v_i)=d_i$ for $1\geqslant i \geqslant n$. This implies that $G$ has order $n$ and size $m$.

  Now, we will show that $G$ is a divisor digraph. Suppose that $( v_i, v_j )$
  and $( v_j, v_k )$ are two edges of $G$. Then from the above construction,
  $1 \leqslant i, j \leqslant d$ and $k \leqslant d_j + 1 \leqslant d_i + 1$.
  Thus $( v_i, v_k )$ is also an edge of $G$. This implies that $G$ is a
  divisor digraph. Let $H$ be the graph obtained from $G$ by ignoring the direction of edges
  of $G$. Then $H$ is a divisor graph of order $n$ and size $m$. This
  concludes the proof of the lemma.
\end{proof}

\begin{lemma}\label{lemma 2}
  Let $n-1 = d_1 \geqslant \ldots \geqslant d_n \geqslant 1$ be a sequence of natural numbers
  and set $e_i = | \{ j : d_j \geqslant i \} |$. Suppose that
  \begin{equation}\label{3} \sum_{i = 1}^n d_i = 2 m < n(n-1), \end{equation}
  and
  \begin{equation}\label{4} \sum_{i = 1}^t d_i = \sum_{i = 1}^t ( e_i - 1 ) \end{equation}
  for all $1 \leqslant t \leqslant f = | \{ i : d_i \geqslant i \} |$. Then
  there exists a sequence $n-1 = d_1' \geqslant \ldots \geqslant d_n' \geqslant 1$ of natural
  numbers such that
  \begin{equation}\label{5} \sum_{i = 1}^n d_i' = 2 ( m + 1 ), \end{equation}
  and
  \begin{equation}\label{6} \sum_{i = 1}^t d_i' = \sum_{i = 1}^t ( e_i' - 1 ) \end{equation}
  for all $1 \leqslant t \leqslant f' = | \{ i : d_i' \geqslant i \} |$, where
  $e_i' = | \{ j : d_j \geqslant i \} |$ for $1 \leqslant i \leqslant n$.
\end{lemma}

\begin{proof}
If $f \geqslant n - 1$ then $n - 1 \geqslant d_1 \geqslant \ldots \geqslant
d_{n - 1} \geqslant d_f \geqslant f = n - 1$. Hence $d_1 = . \ldots = d_{n -
1} = n - 1$, and
\[ \sum_{i = 1}^{n - 1} e_i = n ( n - 1 ) . \]
We have $e_i \leqslant n$ for $1 \leqslant i \leqslant n - 1$, so $e_1 =
\ldots = e_{n - 1} = n$. Hence $d_i \geqslant n - 1$ for $1 \leqslant i
\leqslant n$ or
\[ \sum_{i = 1}^n d_i = n ( n - 1 ), \]
which is a contradiction. Thus, $f < n - 1$. Let $g$ be the smallest index
such that $d_i < n - 1$. Then we have $2 \leqslant g \leqslant f + 1$ (since
$d_{f + 1} \leqslant f < n - 1$). We have two cases.
\begin{enumeratenumeric}
  \item Suppose that $2 \leqslant g \leqslant f$. Set $h = d_g + 2$, $d_g' =
  d_g + 1$, $d_h' = d_h + 1$ and $d_i' = d_i$ for $i \neq g, h$. We have (\ref{5})
  holds since
  \[ \sum_{i = 1}^n d_i' = 2 + \sum_{i = 1}^n d_i = 2 ( m + 1 ) . \]
  Recall that from the proof of Lemma \ref{lemma 1}, we have
  \begin{equation}\label{7}
    d_i = 
    \begin{cases}
    e_i - 1& \tmop{if}\,\, i \leqslant f\\
    f & \tmop{if}\,\, i = f+1\\
    e_{i - 1} &\tmop{if}\,\, i > f + 1.
    \end{cases}
  \end{equation}
  Since $d_g = h - 2 \geqslant g$ and $g \leqslant f$, we have $e_g = h - 1$.
  This implies that $d_{h - 1} \geqslant g$ and $d_h < g$. Besides, $d_1 =
  \ldots = d_{g - 1} = n - 1$ so $d_h \geqslant d_n = e_{n - 1} \geqslant g -
  1$. Hence $d_h = g - 1$, $d_g' = h - 1$ and $d_h' = g$. Therefore, we have
  $e_g' = e_g + 1$, $e_h' = e_h + 1$ and $e_i' = e_i$ for $i \neq g, h$. We
  have $d_h = g - 1 \leqslant f - 1$ so $h > f + 1$. Hence $d_{f + 1}' = d_{f
  + 1} < f + 1$ or $f' \leqslant f$. And we have (\ref{6}) holds for $1\leqslant t \leqslant f'$ since (\ref{4}) holds for $1 \leqslant t \leqslant f' \leqslant f$.
  
  \item Suppose that $g = f + 1$. Then from (\ref{7}) we have $d_g = f = g - 1$. Set
  $h = f + 2$, $d_g' = d_g + 1$, $d_h' = d_h + 1$ and $d_i' = d_i $ for $i
  \neq g, h$. Then it is clear that (\ref{5}) holds. Besides, we have $d_1 = \ldots
  = d_{g - 1} = n - 1$ so
  \[ g - 1 = d_g \geqslant \ldots \geqslant d_n = e_{n - 1} \geqslant g - 1.
  \]
  Hence $d_g = \ldots = d_n = g - 1$. We have $d_g' = d_h' = g < h$, so $f' =
  f + 1$, $e_g' = e_g + 2 = g + 1$ and $e_i' = e_i$ for $i \neq g$. From (\ref{4}),
  we have
  \[ \sum_{i = 1}^t d_i' = \sum_{i = 1}^t ( e_i' - 1 ) \]
  for $1 \leqslant t \leqslant f$. We only need to check for $t = f' ( = f + 1
  = g )$. We have
  \[ \sum_{i = 1}^g d_i' = g + \sum_{i = 1}^f d_i' = e_g' - 1 + \sum_{i = 1}^f
     ( e_i' - 1 ) = \sum_{i = 1}^g ( e_i' - 1 ) . \]
  Thus, (\ref{6}) holds for $1 \leqslant t \leqslant f'$.
\end{enumeratenumeric}
This concludes the proof of the lemma.  
\end{proof}

We are now ready to prove Theorem 1. From Lemma \ref{lemma 1}, we start with the sequence $(n-1,1,\ldots,1)$ to obtain a divisor graph of size $n-1$. Then apply Lemma \ref{lemma 1} and Lemma~\ref{lemma 2} inductively to obtain divisor graphs of order $n,\ldots,\binom{n}{2}$. To construct a divisor graph of order $n$ and size $m$ with $m < n-1$, we choose a vertex and join it with $m$ other vertices. Thus, there exists a divisor graph of order $n$ and size $m$ for any $n$ and $0 \leqslant m \leqslant \binom{n}{2}$. This concludes the proof of the theorem.

\begin{remark}
  An interesting and open question is to find necessary and sufficient
  conditions for a non-increasing sequence $n - 1 \geqslant d_1 \geqslant
  \ldots \geqslant d_n \geqslant 1$ such that there exists a divisor graphs
  with degree sequence $( d_1, \ldots, d_n )$.
\end{remark}

\section{Proof of Theorem \ref{theorem 1}}

  Suppose that $G$ is a divisor graph. Then there exists a set $V$ of positive
  integer such that $G \simeq G ( V )$. We give an orientation on each edge $(
  i, j )$ of $G$ as follows
  \[ ( i, j ) \in E ( G ),\,\, i \rightarrow j \,\,\,\tmop{if} \tmop{and} \tmop{only}
     \tmop{if} \,\,\,i \mid j. \]
  Suppose that $( x, y ), ( y, z )$ are edges of $D$. Then $x \mid y$ and $y \mid z$.
  Hence $x \mid z$ and $( x, z )$ is an edge of $G$.
  
  Now suppose that there exists an orientation $D$ of $G$ such that if $( x,
  y ), ( y, z )$ are edges of $D$ then so is $( x, z )$. We will show that $G$
  is a divisor graph. We will give an explicit labelling for $G$. We start
  with any vertex of $G$ and label it by $\{ a_1 \}$ (a list of one symbol).
  Suppose that we have labelled $k$ vertices of $G$ and we have used $a_1,
  \ldots, a_l$ symbols (each vertex is labelled by a list of symbols and we
  will update this list in each step). We choose any unlabelled vertex, says
  $v$. Consider two sets
  \begin{eqnarray*}
    D_I ( v ) = \{ u \in V ( G ) \mid ( u, v ) \in E ( D ) \}, &  & \\
    D_O ( v ) = \{ u \in V ( G ) \mid ( v, u ) \in E ( D ) \} . &  & 
  \end{eqnarray*}
  We label $v$ by $L_v = \{ a_{l + 1} \}$. For each $u \in D_I ( v )$ and $u$
  was labelled by a list $L_u$ then we add this list into the list $L_v$ to
  have a new list $L_v$ for $v$. And for each $u \in D_O ( v )$ which was
  labelled by a list $L_u$ then we add the new list $L_v$ into $L_u$ to have a
  new list for $u$. For each updated vertex $u$, we consider the set $D_O ( u
  )$. For each $w \in D_O ( u )$ which was labelled  by a list $L_w$, we add
  the new list $L_u$ into the list $L_w$ to have a new list for $w$. We keep
  doing until we have no vertex to update or we come back to some vertex which
  we met along the way. But in the latter case, we have a sequence of
  vertices, says $w_1, \ldots, w_t$ such that $w_1 \rightarrow w_2 \rightarrow
  \ldots \rightarrow w_t \rightarrow w_1$ in $D$. This implies that
  \[ w_1 \mid w_2, w_2 \mid w_3, \ldots, w_t \mid w_1 . \]
  Hence $w_1 = \ldots = w_t$, which is a contradiction. Thus the process must
  be stopped. We repeat the process until all the vertices of $G$ have been
  labelled by lists of symbols. Suppose that we have used $r$ symbols $a_1,
  \ldots, a_r$. We choose $r$ distinct primes $p_1, \ldots, p_r$ and for each
  vertex $v \in V ( G )$ which is labelled by a list $\{ a_{i_1}, \ldots,
  a_{i_v} \} \subseteq \{ a_1, \ldots, a_r \}$ then we label it by the number
  \[ n ( v ) = p_{i_1} \ldots p_{i_v} . \]
  From the construction above, if $( u, v )$ is an edge of $G$ then either
  $L_u \subset L_v$ or $L_v \subset L_u$. This implies that either $u \mid v$ or
  $v \mid u$. Hence $G$ is a divisor graph. This concludes the proof of the
  theorem.
\section{Acknowlegement}
I would like to thank Professor Ping Zhang for sending me reference papers \cite{mot, ba}.

\end{document}